\documentclass[10pt]{article}
\usepackage[english]{babel}
\usepackage[T1]{fontenc}
\usepackage{latexsym}
\usepackage{amssymb}
\usepackage{stmaryrd}
\usepackage{amsmath}
\usepackage{paralist} 
\usepackage[dvips]{graphicx,epsfig}
\usepackage{color}

\newtheorem{thm}{Theorem}

\newtheorem{rem}{Remark}

\newcommand{\conv}[1]{\textrm{conv}}
\newcommand{\var}[1]{\textrm{Var}}
\newcommand{\CQFD}{\nolinebreak \hfill\rule{2mm}{2mm}\medbreak\par}                                                                                

\begin{document}
\begin{center}
{\Large \textbf{On convex hull of $d$-dimensional fractional Brownian motion}}\\
\vspace{0.5cm}

{\large Yu. Davydov\footnote{University  Lille 1, Laboratory P. Painlevé, France}
}
\end{center}

\vspace{1.5cm}

\noindent\textbf{Abstract:} {\small
\vspace{5pt}

It is well known that  for standard Brownian motion $ \{B(t), \;t \geq 0\}$ with values in $\mathbb{R}^d$
its convex hull $ V(t)=\conv \{\{\,B(s),\;s \leq t \}$ with probability $1$ contains $0$ as  an interior point for each $t > 0$ (see \cite{E}). The aim of this note is to state the analoguos property for $d$-dimensional fractional Brownian motion.
 
}

\vspace{10pt}

\vspace{5pt}

\section{Introduction}

Let $(\Omega, \cal{F}, \mathbb{P})$ be a basic probability space. Consider a $d$-dimensional centered Gaussian process $X = \{X(t),\; t\geq 0\}$ defined on $\Omega$ which is self-similar of index $H>0.$
It means that for each constant $c>0$ the process\\ $\{X(ct),\; t\geq 0\}$ has the same distribution as 
$\{c^HX(t),\; t\geq 0\}.$ 
\vspace{5pt}

We call $X$ {\it fractional Brownian motion} (FBM) if for each $e \in \mathbb{R}^d$ the scalar process 
$t \rightarrow \langle X(t),e\rangle$ is a standard one-dimensional FBM up to a constant $c(e).$

It is easy to see that in this case $c^2(e) = \langle Qe,e\rangle,$ where $Q$ is the covariance matrix of $X(1),$ and hence
$$
\mathbb{E}\langle X(t),e\rangle\langle X(s),e\rangle = \langle Qe,e\rangle\frac{1}{2}(t^{2H} +s^{2H}-|t-s|^{2H}),\;\;\; t,s \geq 0;\,\, e\in \mathbb{R}^d
$$
(see \cite{RS},\cite{LPS} and references therein for more general definitions of operator self-similar FBM).

The next  properties follow from the definition without difficulties.

\begin{enumerate}
\item[1)] {\bf Continuity.} The process $X$ has a continuous version.

Below we always suppose $X$ to be continuous.

\item[2)] {\bf Reversibility.} Define the process $Y$ by
$$
Y(t) = X(1) -X(1-t),\;\;\; t\in [0,1].
$$
Then $\{Y(t),\; t\in [0,1]\} \stackrel{\cal L}{=} \{X(t),\; t\in [0,1]\},$ 
where $\stackrel{\cal L}{=}$ means equality in law.

\item[3)] {\bf Ergodicity.} Let $L = \{L(u),\; u\in \mathbb{R}^1\}$ be the strictly stationary Gaussian process obtained from $X$ by Lamperti transformation:
\begin{equation}
L(u) = e^{-Hu} X(e^u),\;\;\; u \in \mathbb{R}^1.
\end{equation}\label{lamperti}
Then $L$ is ergodic (see \cite{CFS}, Ch. 14, \S 2,  Th.1, Th.2).
\end{enumerate}

\section{Results}

For Borel set $A \subset \mathbb{R}^d$ we denote by $\mathrm{conv}{(A)} $ the closed convex hull generated by $A$.

The object of our interest is the convex hull process related to $X$:
$$
V(t) = \mathrm{conv} \{\,X(s),\;\;s\leq t\}.
$$

It is supposed below  that the law of $X(1)$ is non degenerate, that is the rank of the matrix $Q$ is equal to $d$.

\begin{thm} \label{th1} 
With probability $1$ for all $t>0$ the point $0$ is an interior point of $V(t).$
\end{thm}

As a corollary we immediately deduce the following fact.

\begin{thm} \label{th2} 
For each $t>0$ with probability $1$  the point $X(t)$ is an interior point of $V(t).$
\end{thm}

{\bf Proof of Th. 2.} Denote by $A^\circ $ the interior of $A.$ By self-similarity of the process $X$ it is sufficient to state this property for $t=1.$ Then, due to the reversibility of $X$ by Th. 1., a.s. 
\begin{equation}
\label{revers}
0 \in [\mathrm{conv}\{\,X(1) -X(1-t),\;\;\; t\in [0,1]\}]^\circ.
\end{equation}

As
$$
\mathrm{conv}\{\,X(1) -X(1-t),\;\;\; t\in [0,1]\} = X(1) -\conv\{\{\,X(s),\;\;\; s\in [0,1]\},
$$
the relation  (\ref{revers}) is equivalent to
$$
X(1) \in [\mathrm{conv}\{\,X(s),\;\;s\in [0,1]\}]^\circ,
$$
which concludes the proof.\CQFD

Let ${\cal K}_d$ be the family of all compact convex subsets of $\mathbb{R}^d.$ It is well known that
${\cal K}_d$ equipped with Hausdorff metric is a Polish space.
\vspace{5pt}

We say that a function $f : [0,1] \rightarrow {\cal K}_d$ is {\it increasing}, if $f(t) \subset f(s)$ for
$0\leq t<s\leq 1.$
\vspace{5pt}

We say that a function $f : [0,1] \rightarrow {\cal K}_d$ is {\it Cantor - staircase} (C-S), 
if $f$ is continuous, increasing and such that for almost every $t\in [0,1]$ there exists an interval $(t-\varepsilon, t+\varepsilon)$ where $f$ is constant.
\vspace{5pt}

The next statement is an easy corollary of Th.2.

\begin{thm} \label{th3} 
With probability $1$ the paths of the process $t\rightarrow V(t)$ are C-S functions.

\end{thm}

{\bf Proof of Th. 3.} We use the notation $X(s,\omega)$ for $X(s)$  to emphasize the dependance of $\omega \in \Omega.$ By Th.2 for each $t\in (0,1)$ with probability 1 $X(t) \in V(t)^\circ.$
By continuity of $X$, for $\mathbb{P}$-almost each $\omega$ there exists $\varepsilon >0$ such that
$X(s,\omega) \in V(t)^\circ$ for all $s \in (t-\varepsilon, t+\varepsilon)$ which gives
$V(s) = V(t),\;\;\;$\\$ \forall s \in (t-\varepsilon, t+\varepsilon).$
\CQFD

\begin{rem} \label{th1-2}
Let $h :{\cal K} \rightarrow \mathbb{R}^1 $ be an increasing continuous function. Then  almost all paths of the process $t \rightarrow h(V(t))$ are C-S real  functions. This obvious fact may be applied to all reasonable geometrical caracteristics of $V(t),$ such as volume, surface area, diameter,...
\end{rem}

\section{Proof of Theorem 1} Let $\Theta = \{0,1\}^d$ be the set of all diadic sequences of length $d.$
Denote by $D_\theta,$\\$ \;\theta \in \Theta,$ the quadrant 
$$
D_\theta = \prod_{i=1}^{d}\mathbb{R}_{\theta_i},
$$
where $\mathbb{R}_{\theta_i} = [0,\infty)$ if  $\theta_i = 1,$ and $\mathbb{R}_{\theta_i} = (-\infty, 0]$ if
$\theta_i = 0.$

The positive quadrant $D_{(1,1,\ldots,1)}$ for simplicity is denoted by $D.$

\noindent
 We  first show that
\begin{equation}
p\;\; \stackrel{def}{=}\;\; \mathbb{P}\{\,\exists\, t\in (0,1]\;\vert\; X(t) \in D^\circ\} = 1.
\end{equation}
Remark that $p$ is strictly positive:
\begin{equation}
\label{p}
p\geq \mathbb{P}\{X(1) \in D^\circ\} > 0
\end{equation}
due to the hypothesis that the law of $X(1)$ is non degenerate.

By self similarity
$$
\mathbb{P}\{D^\circ\ \cap \{X(t), t \in [0,T]\}= \emptyset\} = 1-p
$$
for every $T>0.$

The sequence of events $(A_n)_{n\in\mathbb{N}}\,,$
$$
A_n = \{D^\circ\ \cap \{X(t), t \in [0,n]\}= \emptyset\},
$$
being decreasing, it follows that
$$
1-p = \lim \mathbb{P}(A_n) = \mathbb{P}(\cap_nA_n) = \mathbb{P}\{X(t) \notin D^\circ,\;\; \forall t\geq0\}.
$$
In terms of the stationary process $L$ from Lamperti representation (\ref{lamperti}) it means that
$$
\mathbb{P}\{L(s) \notin D^\circ,\;\; \forall s \in \mathbb{R}^1 \}= 1-p.
$$
As this event is invariant, by ergodicity of $L$ and due to (\ref{p}) we see that the value $p=1$ is the  only one possible.

Applying the analoguos arguments to another quadrants $D_\theta$, we get that with probability 1 there exists points  $t_\theta \in (0,1],$ such that\\ $X(t_\theta) \in D_\theta^\circ, \;\;\;\forall\theta \in \Theta.$
Now, to end the proof it is sufficient to remark that
$$
V(1)^\circ = \mathrm{conv}\{X(t), t \in [0,1]\}^\circ \supset \mathrm{conv} \{X(t_\theta), \theta \in \Theta\}^\circ
$$
and that the last set evidently contains $0.$ \CQFD
\vspace{10pt}

{\bf Acknowledgments.} The author wishes to thank 
 all participants of working seminar on stochastic geometry of the university Lille 1 for their support.


\end{document}